\documentclass[12pt]{amsart}
\usepackage{amssymb,amsmath}
\def\eps{\varepsilon}
\def\e{{\rm e}}
\def\d{{\rm d}}
\def\ddt{\frac{\d}{\d t}}

\def\R {\mathbb{R}}
\def\N {\mathbb{N}}
\def\BB {\mathfrak{B}}
\def\H {{\mathcal H}}
\def\V {{\mathcal V}}

\def\C {{\mathcal C}}
\def\D {{\mathcal D}}

\def\G {{\mathcal G}}
\def\E {{\mathcal E}}
\def\L {{\mathcal L}}
\def\A {{\mathfrak A}}

\def\J {{\mathcal J}}

\def\Q {{\mathcal Q}}

\def\W {{\mathcal W}}
\def \l {\langle}
\def \r {\rangle}
\def \pt {\partial_t}
\def \ptt {\partial_{tt}}

\newtheorem{proposition}{Proposition}[section]
\newtheorem{theorem}[proposition]{Theorem}
\newtheorem{corollary}[proposition]{Corollary}
\newtheorem{lemma}[proposition]{Lemma}
\theoremstyle{definition}

\newtheorem{remark}[proposition]{Remark}

\numberwithin{equation}{section}

\def \au {\rm}
\def \ti {\it}
\def \jou {\rm}
\def \bk {\it}
\def \no#1#2#3 {{\bf #1} (#3), #2.}
\def \eds#1#2#3 {#1, #2, #3.}

\title[Extensible Thermoelastic Beam System]
{Global Attractors for the Extensible \\
Thermoelastic Beam System}

\author[C. Giorgi, M.G. Naso, V. Pata, M. Potomkin]
{C. Giorgi, M.G. Naso, V. Pata, M. Potomkin}

\address{Universit\`a di Brescia - Dipartimento di Matematica
\newline\indent
Via Valotti 9, 25133 Brescia, Italy}
\email{giorgi@ing.unibs.it}
\email{naso@ing.unibs.it}

\address{Politecnico di Milano - Dipartimento di Matematica ``F.\ Brioschi''
\newline\indent
Via Bonardi 9, 20133 Milano, Italy}
\email{vittorino.pata@polimi.it}

\address{Kharkov National University - Department of Mathematics and Mechanics
\newline\indent
4 Svobody sq, 61077 Kharkov, Ukraine}
\email{mika\_potemkin@mail.ru}

\thanks{Work partially supported by the Italian PRIN research project 2006
{\it Problemi a frontiera libera, transizioni di fase e modelli
di isteresi}.}



\begin{document}

\begin{abstract}
This work is focused on the dissipative system
$$
\begin{cases}
\ptt u+\partial_{xxxx}u
+\partial_{xx}\theta-\big(\beta+\|\partial_x u\|_{L^2(0,1)}^2\big)\partial_{xx}u=f\\
\noalign{\vskip.7mm}
\pt \theta -\partial_{xx}\theta -\partial_{xxt} u= g
\end{cases}
$$
describing the dynamics
of an extensible thermoelastic beam, where
the dissipation is entirely contributed
by the second equation ruling the evolution of $\theta$.
Under natural boundary conditions, we prove the existence
of bounded absorbing sets.
When the external sources $f$ and $g$ are time-independent,
the related semigroup of solutions is shown to possess the global
attractor of optimal regularity for all parameters $\beta\in\R$.
The same result holds true when the first equation
is replaced by
$$
\ptt u-\gamma\partial_{xxtt} u+\partial_{xxxx}u
+\partial_{xx}\theta-\big(\beta+\|\partial_x u\|_{L^2(0,1)}^2\big)\partial_{xx}u=f
$$
with $\gamma>0$. In both cases, the solutions on the attractor are strong solutions.
\end{abstract}

\maketitle

\section{Introduction}

\noindent
For $t>0$, we consider the evolution system
\begin{equation}
\label{PROB}
\begin{cases}
\displaystyle
\ptt u+\partial_{xxxx}u
 +\partial_{xx}\theta-\Big(\beta+\int_0^1
|\partial_x u(x,\cdot)|^2\d x\Big)\partial_{xx}u=f,\\
\noalign{\vskip.7mm}
\pt \theta -\partial_{xx}\theta -\partial_{xxt} u= g,
\end{cases}
\end{equation}
in the unknown variables $u=u(x,t):[0,1]\times\R^+\to\R$
and $\theta=\theta(x,t):[0,1]\times\R^+\to\R$, having
put $\R^+=[0,\infty)$.
The two equations are supplemented
with the initial conditions
\begin{equation}
\label{IC}
\begin{cases}
u(x,0)=u_0(x),\\
\pt u(x,0)=u_1(x),\\
\theta(x,0)=\theta_0(x),
\end{cases}
\end{equation}
for every $x\in[0,1]$, where $u_0$, $u_1$ and $\theta_0$ are assigned data.

System~\eqref{PROB} describes
the vibrations of an extensible thermoelastic beam
of unitary natural length, and is obtained by combining
the pioneering ideas of Woinowsky-Krieger \cite{W} with the theory
of linear thermoelasticity \cite{CAR}.
Although a rigorous variational derivation of the model will 
be addressed in a forthcoming paper,
it is worth noting that \eqref{PROB} is a 
mild quasilinear version of the nonlinear motion equations devised in 
\cite[\S3]{LLS}.

With regard to the physical meaning of the variables
in play, $u$ represents the vertical deflection of the beam from
its configuration at rest, while the ``temperature variation"
$\theta$ actually arises from an approximation of the temperature variation with 
respect to a reference value, and it has the dimension of a 
temperature gradient (see \cite{LLS}).
The real function $f=f(x,t)$ is the lateral load distribution
and $g=g(x,t)$ is the external heat supply, having the role
of a control function.
Finally, the parameter $\beta\in\R$ accounts for the axial force acting in
the reference configuration:
$\beta>0$ when the beam is stretched, $\beta<0$
when compressed.
Concerning the boundary conditions, for all $t\geq 0$ we assume
\begin{equation}
\label{BC}
\begin{cases}
u(0,t)=u(1,t)=\partial_{xx}u(0,t)=\partial_{xx}u(1,t)=0,\\
\theta(0,t)=\theta(1,t)=0.
\end{cases}
\end{equation}
Namely, we take Dirichlet boundary conditions for the temperature variation $\theta$
and hinged boundary conditions for the vertical
deflection $u$.

The focus of this paper is the study of the longterm properties
of the dynamical system generated by problem \eqref{PROB}-\eqref{BC}
in the natural weak energy phase space.
In particular, for the autonomous case, we prove
the existence of the global attractor of optimal regularity, for all values
of the real parameter $\beta$.
The main difficulty arising in the asymptotic analysis
comes from the very weak dissipation exhibited by the model,
entirely contributed by the thermal component,
whereas the mechanical component, by itself, does not
cause any decrease of energy.
Hence, the loss of mechanical energy is due
only to the coupling, which propagates the thermal dissipation
to the mechanical component with the effect of producing
mechanical dissipation. From the mathematical side,
in order to obtain
stabilization properties it is necessary
to introduce sharp energy functionals,
which allow to exploit the thermal dissipation in its full strength.
A similar situation has been faced in
\cite{AL1,AL2,AL3,CL,CHLA02,CHLA08,GIOP,LAS,POT}, dealing with
linear and semilinear thermoelastic problems
without mechanical dissipation. 
Along the same line, we also mention
\cite{LMS}, which considers a quasilinear thermoelastic plate system.

After this work was finished, we learned of a paper by Bucci and Chueshov~\cite{BUC},
which treats (actually, in a more general version) the same problem
discussed here. In~\cite{BUC}, borrowing some techniques from the recent article \cite{CHLA08},
the authors prove the existence
of the global attractor of optimal regularity and finite fractal dimension
for the semigroup generated by the autonomous version
of \eqref{PROB}-\eqref{BC}. The existence of the (regular) attractor is also shown
in presence of a rotational term in the first equation, whose dynamics
becomes in turn of hyperbolic type. 
The proofs of~\cite{BUC} heavily rely on two basic facts: a key estimate,
nowadays known in the literature as {\it stabilizability inequality}
(cf.\ \cite{CL,CHLA08}), and
the gradient system structure featured by the model. The regularity of the attractor
is demonstrated only in a second moment,
exploiting the peculiar form
of the attractor itself: a section of all bounded complete orbits of the semigroup.

Nevertheless, we still believe that our paper might be of some interest,
at least for the following 
reasons:
\begin{itemize}
\item[(i)] We do not appeal to
the gradient system structure, except for the characterization
of the attractor as the unstable set of stationary solutions. 
Accordingly, the existence of absorbing sets is established via
explicit energy estimates, providing a precise (uniform) control
on the entering times of the trajectories. The method applies to
the nonautonomous case as well, where the gradient system structure is lost.
\item[(ii)] Our proof of asymptotic compactness is rather direct and simpler
than in \cite{BUC}. Indeed, it merely boils down to the construction of
a suitable decomposition of the semigroup. Incidentally, the required regularity is gained
in just one single step, without making use of bootstrapping arguments.
\item[(iii)] As a matter of fact, we prove a stronger result:
we find exponentially
attracting sets of optimal regularity, obtaining at the same time the
attractor and its regularity.
Having such
exponentially attracting sets,
it is possible to show with little effort the existence of regular
exponential attractors in the sense of \cite{EMZ}, having
finite fractal dimension.
\item[(iv)] In the rotational case, where the first equation
contains the extra term
$-\gamma\partial_{xxtt}u$ with $\gamma>0$,
we improve the regularity of the attractor devised in \cite{BUC}.
\end{itemize}

\subsection*{Plan of the paper}
In Section~2, we consider an abstract generalization
of \eqref{PROB}-\eqref{BC}, whose solutions are generated
by a family of solution operators $S(t)$.
We also discuss the changes needed in the abstract framework
to take into account other types of boundary conditions
for $u$.
Section~3 is devoted to the existence
of an absorbing set for $S(t)$. In Section~4, we dwell on the autonomous case,
where $S(t)$ is a semigroup, establishing the existence 
and the regularity of the
global attractor.
The proofs
are carried out in the subsequent  Section~5.
In the last Section~6, we extend the results to a more general model,
where an additional rotational inertia term is present.

\section{The Abstract Problem}

\subsection{Notation}
Let $(H,\l\cdot,\cdot\r,\|\cdot\|)$ be a separable real Hilbert space,
and let $A:H\to H$ be a strictly positive selfadjoint operator
with domain $\D(A)\Subset H$.
For $r\in\R$, we introduce the scale of Hilbert spaces generated by the powers of
$A$
$$
H^r=\D(A^{r/4}),\quad
\l u,v\r_r=\l A^{r/4}u,A^{r/4}v\r,\quad
\|u\|_r=\|A^{r/4}u\|.
$$
We will always omit the index $r$ when $r=0$. The symbol $\l\cdot,\cdot\r$ will
also be used to denote the duality product between $H^r$ and
its dual space $H^{-r}$.
In particular, we have the compact embeddings $H^{r+1}\Subset H^r$, along with
the generalized Poincar\'e inequalities
$$\lambda_1\|u\|_r^4\leq \|u\|_{r+1}^4,\quad\forall u\in H^{r+1},$$
where $\lambda_1>0$ is the first eigenvalue of $A$. Finally, we define
the the product Hilbert spaces
$$\H^r=H^{r+2}\times H^r\times H^r.$$

\subsection{Formulation of the problem}
For $\beta\in\R$, we consider the abstract Cauchy problem on
$\H$
in the unknown variables $u=u(t)$ and $\theta=\theta(t)$
\begin{equation}
\label{BASE}
\begin{cases}
\ptt u+Au-A^{1/2}\theta+
\big(\beta+\|u\|^2_1\big)A^{1/2}u= f(t),\quad t>0,\\
\noalign{\vskip.7mm}
\pt \theta+A^{1/2}\theta+A^{1/2}\pt u=g(t),\quad t>0,\\
u(0)=u_0,\quad
\pt u(0)=u_1,
\quad\theta(0)=\theta_0.
\end{cases}
\end{equation}
The following well-posedness result holds.

\begin{proposition}
\label{EU}
Assume that
$$f\in L^1_{\rm loc}(\R^+,H),\quad
g\in L^1_{\rm loc}(\R^+,H)+L^2_{\rm loc}(\R^+,H^{-1}).$$
Then, for all initial data $(u_0,u_1,\theta_0)\in\H$, problem~\eqref{BASE}
admits a unique (weak) solution
$$(u(t),\pt u(t),\theta(t))\in\C(\R^+,\H)$$
with
$$(u(0),\pt u(0),\theta(0))=(u_0,u_1,\theta_0).$$
Moreover, calling $\bar z(t)$ the difference of 
any two solutions
corresponding to initial data having norm less than or equal to $R\geq 0$,
there exists $C=C(R)\geq 0$ such that
$$
\|\bar z(t)\|_\H
\leq C\e^{C t}\|\bar z(0)\|_\H,\quad\forall t\geq 0.
$$
\end{proposition}

We omit the proof, based on
a standard Galerkin approximation procedure together with a slight generalization
of the usual Gronwall lemma (cf.\ \cite{PPV}).
Proposition~\ref{EU} translates into the existence of
the {\it solution operators}
$$S(t):\H\to\H,\quad t\geq 0,$$
acting as
$$z=(u_0,u_1,\theta_0)\mapsto S(t)z=(u(t),\pt u(t),\theta(t)),$$
and satisfying the joint continuity property
$$(t,z)\mapsto S(t)z\in\C(\R^+\times\H,\H).$$

\begin{remark}
In the autonomous case, namely, when both $f$ and $g$ are time-independent,
the family $S(t)$ fulfills the semigroup property
$$S(t+\tau)=S(t)S(\tau),\quad\forall t,\tau\geq 0.$$
Thus,
$S(t)$ is a strongly continuous semigroup of operators on $\H$.
\end{remark}

We define the {\it energy} at time $t\geq 0$ corresponding
to the initial data $z=(u_0,u_1,\theta_0)\in\H$ as
$$\E(t)=\frac12\|S(t)z\|^2_\H+\frac14\big(\beta+\|u(t)\|_1^2\big)^2.
$$
Multiplying the first equation
of \eqref{BASE} by $\pt u$ and the
second one by $\theta$, we find the {\it energy identity}
\begin{equation}
\label{E}
\ddt \E+\|\theta\|_1^2=\l \pt u,f\r+\l\theta,g\r.
\end{equation}
Indeed,
\begin{align*}\frac14\ddt\big(\beta+\|u\|_1^2\big)^2
&=\frac12\big(\beta+\|u\|_1^2\big)
\ddt\big(\beta+\|u\|_1^2\big)\\
&=\frac12\big(\beta+\|u\|_1^2\big)\ddt \|u\|_1^2
=\big(\beta+\|u(t)\|_1^2\big)\l A^{1/2}u,\pt u\r.
\end{align*}
As a consequence, for every $T>0$, there exists a positive increasing
function $\Q_T$ such that
\begin{equation}
\label{TFIN}
\E(t)\leq \Q_T(\E(0)),\quad\forall t\in [0,T].
\end{equation}

\subsection{The concrete problem: other boundary conditions}

The abstract system~\eqref{BASE}
serves as a model to describe quite general situations, including
thermoelastic plates. 
In particular, problem~\eqref{PROB}-\eqref{BC} 
is a concrete realization of \eqref{BASE},
obtained by putting $H=L^2(0,1)$ and $A=\partial_{xxxx}$ with domain
$$\D(A)=\big\{u\in H^4(0,1):u(0)=u(1)=\partial_{xx} u(0)=\partial_{xx} u(1)=0\big\}.
$$
In which case,
$$A^{1/2}=-\partial_{xx},\quad\D(A^{1/2})
=H^2(0,1)\cap H_0^1(0,1).
$$
However,
although we supposed in \eqref{BC} that both ends of the beam are hinged,
different boundary conditions
for $u$ are physically significant as well, such as
\begin{equation}
\label{CLAMP}
u(0,t)=u(1,t)=\partial_{x}u(0,t)=\partial_{x}u(1,t)=0,
\end{equation}
when both ends of the beam are clamped, or
\begin{equation}
\label{MIX}
u(0,t)=u(1,t)=\partial_{x}u(0,t)=\partial_{xx}u(1,t)=0,
\end{equation}
when one end is clamped and the other one is hinged.
On the contrary, the so-called cantilever boundary condition 
(one end clamped and the other one free) 
does not comply with the extensibility assumption, since 
no geometric constraints compel the beam length to change.

In order to write an abstract formulation accounting also
for the boundary conditions \eqref{CLAMP} and \eqref{MIX},
let $A_\star:\D(A_\star)\Subset H\to H$
be another selfadjoint strictly positive operator.
Accordingly,
for $r\in\R$, we have the further scale of Hilbert
spaces
$$
H^r_\star=\D(A_\star^{r/4}),\quad
\l u,v\r_{r,\star}=\l A^{r/4}_\star u,A^{r/4}_\star v\r,\quad
\|u\|_{r,\star}=\|A^{r/4}_\star u\|.
$$
Defining
$$\H^r_\star=H^{r+2}_\star\times H^r\times H^r,$$
we consider the evolution system
\begin{equation}
\label{BASEOTHER}
\begin{cases}
\ptt u+A_\star u-A^{1/2}\theta+
\big(\beta+\|u\|^2_1\big)A^{1/2}u= f,\\
\noalign{\vskip.7mm}
\pt \theta+A^{1/2}\theta+A^{1/2}\pt u=g,\\
\end{cases}
\end{equation}
with initial data $(u_0,u_1,\theta_0)\in\H_\star$, and we look for solutions
$$(u(t),\pt u(t),\theta(t))\in\C(\R^+,\H_\star).$$
If there exists
a dense subspace $D$ of $H$, contained in $\D(A_\star)\cap\D(A)$ and such that
\begin{equation}
\label{relabell}
A_\star u=A u\in D,\quad\forall u\in D,
\end{equation}
we can work within a suitable Galerkin approximation scheme, 
using an orthonormal basis
of $H$ made of elements in $D$. Then,
exploiting the equality
$$\|u\|^2_{2,\star}=\l A_\star u,u\r=\l A u,u\r=\|u\|^2_{2},$$
and, in turn, the interpolation inequalities
$$
\|u\|^2_{1}
\leq \eps\|u\|_{2,\star}^2+\frac1{4\eps}\|u\|^2,
\quad
\|u\|^2_{3}
\leq \eps\|u\|_{4,\star}^2+\frac1{4\eps}\|u\|_{2,\star}^2,
$$
valid for every $u\in D$,
we can adapt the proofs
of the subsequent sections,
in order to establish to the existence of the global attractor.

This abstract scheme applies to the concrete case of the extensible
thermoelastic beam with boundary conditions for $u$ of the form
\eqref{CLAMP} or \eqref{MIX}, upon choosing $H=L^2(0,1)$,
$A_\star=\partial_{xxxx}$ with domain
$$\D(A_\star)
=\begin{cases}
H^4(0,1)\cap H^2_0(0,1) & \text{b.c.\ \eqref{CLAMP}},\\
\big\{u\in H^4(0,1)\cap H_0^1(0,1):\partial_{x}u(0)=\partial_{xx}u(1)=0\big\}
& \text{b.c.\ \eqref{MIX}},
\end{cases}
$$
and setting, for instance, 
$D=\C^\infty_{{\rm cpt}}(0,1)$.

\section{The Absorbing Set}

\noindent
In this section, we prove the existence of
an absorbing set for the family $S(t)$.
This is a bounded set $\BB\subset \H$ with the following property:
for every $R\geq 0$, there is an {\it entering time} $t_R\geq 0$ such that
$$\bigcup_{t\geq t_R}S(t)z\subset \BB,
$$
whenever $\|z\|_\H\leq R$.

In fact, we establish a more general result.

\begin{theorem}
\label{thmABS}
Let $f\in L^\infty(\R^+,H)$, and let $\pt f$ and $g$ be translation
bounded functions in $L^2_{\rm loc}(\R^+,H^{-1})$, that is,
\begin{equation}
\label{tb}
\sup_{t\geq 0}\int_t^{t+1}\big\{\|\pt f(\tau)\|^2_{-1}
+\|g(\tau)\|^2_{-1}\big\}\d\tau<\infty.
\end{equation}
Then, there exists $R_0>0$ with the following property:
in correspondence of every $R\geq 0$, there is $t_0=t_0(R)\geq 0$ such that
$$\E(t)\leq R_0,\quad \forall t\geq t_0,
$$
whenever $\E(0)\leq R$. Both $R_0$ and $t_0$ can be explicitly computed.
\end{theorem}

The absorbing set, besides giving
a first rough estimate of the dissipativity of the system,
is the preliminary step to prove the existence
of much more interesting objects describing the asymptotic dynamics,
such as global or
exponential attractors (see, for instance, \cite{BV,CV,CVbook,Chbook,HAL,MZ,TEM}).
Unfortunately, in certain situations where the dissipation
is very weak, a direct proof of the existence of the absorbing set via {\it explicit}
energy estimates
might be very hard to find. On the other hand, for a quite general class of
autonomous problems
(the so-called gradient systems), it is possible
to use an alternative approach and
overcome this obstacle, appealing to the existence of
a Lyapunov functional (see \cite{CP,HAL,LAD}). In which case, if the
semigroup possesses suitable
smoothing properties, one obtains right away the global attractor,
and the absorbing set is then recovered
as a byproduct. Though,
the procedure
provides no quantitative information on the entering time $t_R$,
which is somehow unsatisfactory,
especially in view of numerical simulations.
This technique has been successfully adopted
in the recent paper~\cite{GPV}, concerned with the longterm
analysis of an integrodifferential
equation with low dissipation,
modelling the transversal motion of an extensible viscoelastic beam.

As mentioned in the introduction,
the problem considered in the present work is also weakly dissipative.
But if we assume $f$ and $g$ independent of time,
there is a way to define a Lyapunov functional (actually, for an equivalent problem),
which would allow to exploit the method described above.
In any case, in order to exhibit an actual bound on $t_R$,
and also to deal with time-dependent external forces,
a direct proof of Theorem~\ref{thmABS} would be much more desirable.
However, due to presence of the coupling term,
when performing the standard (and unavoidable)
estimates, some ``pure" energy terms having a power strictly greater than one
pop up with the wrong sign. Such terms cannot be handled
by means of standard Gronwall-type lemmas.
Nonetheless, we are still able to establish the result, leaning on
the following novel Gronwall-type lemma with parameter devised in~\cite{GPZ}.

\begin{lemma}
\label{superl}
Let $\Lambda:\R^+\to\R^+$ be an absolutely continuous function
satisfying, for some $K\geq 0$, $Q\geq 0$, $\eps_0>0$
and every $\eps\in(0,\eps_0]$,
the differential inequality
$$
\ddt\Lambda(t)+\eps \Lambda(t)\leq K\eps^2 [\Lambda(t)]^{3/2}+\eps^{-2/3}\varphi(t),
$$
where $\varphi:\R^+\to\R^+$ is any locally summable function such that
$$\sup_{t\geq 0}\int_t^{t+1}\varphi(\tau)\d \tau\leq Q.$$
Then, there exist $R_1>0$ and $\kappa>0$ such that,
for every $R\geq 0$, it follows that
$$\Lambda(t)\leq R_1,\quad\forall t\geq R^{1/\kappa}(1+\kappa Q)^{-1},$$
whenever $\Lambda(0)\leq R$.
\end{lemma}

\begin{remark}
Both $R_1$ and $\kappa$ can be explicitly computed
in terms of the constants $K,Q$ and $\eps_0$ (cf.\ \cite{GPZ}).
\end{remark}

We are now ready to proceed to the proof of the theorem.

\begin{proof}[Proof of Theorem \ref{thmABS}]
Here and in the sequel,
we will tacitly use several times the Young and the H\"older
inequalities, besides the usual Sobolev embeddings.
The generic positive constant $C$ appearing in this proof
may depend on $\beta$ and $\|f\|_{L^\infty(\R^+,H)}$.

On account of \eqref{E}, the functional
$$\L(t)=\E(t)-\l u(t),f(t)\r
$$
satisfies the differential equality
$$\ddt \L+\|\theta\|_1^2=-\l u,\pt f \r+\l\theta, g\r.$$
Observing that
\begin{equation}
\label{immediate}
\|u\|_1\leq C|\beta+\|u\|_1^2|^{1/2}+C|\beta|^{1/2}\leq C\E^{1/4}+C,
\end{equation}
we have the control
$$-\l u,\pt f\r\leq C\eps^{2/3} \E^{1/2}+\eps^{-2/3}\|\pt f\|^2_{-1}
\leq C\eps^2 \E^{3/2}+\eps^{-2/3}\|\pt f\|^2_{-1}+C,
$$
for all $\eps\in(0,1]$. Moreover,
$$\l \theta,g\r\leq \frac12\|\theta\|_1^2+C\|g\|^2_{-1}.
$$
Thus, we obtain the differential inequality
\begin{equation}
\label{L}
\ddt \L+\frac12\|\theta\|_1^2
\leq C\eps^2 \E^{3/2}+\eps^{-2/3}\|\pt f\|^2_{-1}+C\|g\|^2_{-1}+C.
\end{equation}
Next, we consider the auxiliary functionals
$$\Phi(t)=\l\pt u(t),u(t)\r,\quad
\Psi(t)=\l \pt u(t),\theta(t)\r_{-1}.$$
Concerning $\Phi$, we have
$$
\ddt\Phi+\|u\|_2^2+\big(\beta+\|u\|_1^2\big)^2
-\beta\big(\beta+\|u\|_1^2\big)
=\|\pt u\|^2+\l u,\theta\r_1+\l u,f\r.
$$
Noting that
$$\frac12\big(\beta+\|u\|_1^2\big)^2
-\beta\big(\beta+\|u\|_1^2\big)=\frac12\|u\|_1^4-\frac12 \beta^2\geq -C,$$
and
$$\l u,\theta\r_1+\l u,f\r\leq \frac14\|u\|_2^2+C\|\theta\|_1^2+C,
$$
we are led to
\begin{equation}
\label{PSI1}
\ddt\Phi+\frac34\|u\|_2^2+\frac12\big(\beta+\|u\|_1^2\big)^2
\leq \|\pt u\|^2+C\|\theta\|_1^2+C.
\end{equation}
Turning to $\Psi$, we have the differential equality
$$\ddt\Psi+\|\pt u\|^2=\|\theta\|^2
-\l\pt u,\theta\r-\l u,\theta\r_1+\l\pt u,g\r_{-1}+\l\theta,f\r_{-1}+\J,$$
having put
$$\J=-\big(\beta+\|u\|_1^2\big)\l u,\theta\r.$$
We easily see that
\begin{align*}
&\|\theta\|^2
-\l\pt u,\theta\r-\l u,\theta\r_1+\l g,\pt u\r_{-1}+\l f,\theta\r_{-1}\\
&\leq\frac{1}{8}\|u\|_2^2 + \frac14 \|\pt u \|^2 +C\|\theta\|_1^2+C\|g\|^2_{-1}+C,
\end{align*}
whereas, in light of~\eqref{immediate},
the remaining term $\J$ is controlled as
$$\J\leq C\|\theta\|\big(\|u\|_1^3+1\big)\leq
C\|\theta\|\E^{3/4} +C\|\theta\|
\leq C\|\theta\|\E^{3/4}+C\|\theta\|_1^2+C.$$
In conclusion,
\begin{equation}
\label{PSI2}
\ddt\Psi+\frac34\|\pt u\|^2
\leq \frac{1}{8}\|u\|_2^2+C\|\theta\|_1^2+C\|\theta\|\E^{3/4}+C\|g\|^2_{-1}+C.
\end{equation}
Collecting \eqref{PSI1}-\eqref{PSI2}, we end up with
\begin{equation}
\label{UNODUE}
\ddt\big\{\Phi+2\Psi\big\}
+\E\leq C\|\theta\|_1^2+C\|\theta\|\E^{3/4}+C\|g\|^2_{-1}+C.
\end{equation}
Finally, for $\eps\in(0,1]$, we set
$$\Lambda(t)=\L(t)+2\eps\big\{\Phi(t)+2\Psi(t)\big\}+C,$$
where the above $C$ is large enough and $\eps$ is small
enough such that
\begin{equation}
\label{CTRL}
\frac12 \E\leq\Lambda\leq 2\E+C.
\end{equation}
Then, calling
$$\varphi(t)=C+C\|\pt f(t)\|^2_{-1}+C\|g(t)\|^2_{-1},
$$
the inequalities
\eqref{L}, \eqref{UNODUE} and \eqref{CTRL} entail
\begin{align*}
\ddt\Lambda
+\eps\Lambda +\frac12(1-C\eps)\|\theta\|_1^2
&\leq C\eps^2 \Lambda^{3/2}+C\eps\|\theta\|\Lambda^{3/4}+\eps^{-2/3}\varphi\\
&\leq C\eps^2 \Lambda^{3/2}+\eps^{-2/3}\varphi
+\frac14\|\theta\|_1^2.
\end{align*}
It then is apparent that there exists $\eps_0>0$ small such that, for every $\eps\in(0,\eps_0]$,
$$
\ddt\Lambda
+\eps\Lambda \leq C\eps^2 \Lambda^{3/2}+\eps^{-2/3}\varphi.
$$
By virtue of \eqref{tb}, we are in a position to apply
Lemma~\ref{superl}.
Using once more \eqref{CTRL}, the proof is finished.
\end{proof}

\section{The Global Attractor}

\noindent
In the sequel, we will assume the external forces $f$ and $g$
to be independent of time.
In which case, $S(t)$ is a strongly continuous semigroup on $\H$.
We define
$$\theta_g=A^{-1/2}g,\quad z_g=(0,0,\theta_g).$$
The main result, which will be proved in the next
section, reads as follows.

\begin{theorem}
\label{MAIN}
Let $f,g\in H$.
Then, the semigroup $S(t)$ acting on $\H$ possesses the (connected)
global attractor $\A$. Moreover,
$$\A=z_g+\A_0,$$
where $\A_0$ is a bounded subset of the space
$\H^2\Subset\H$.
\end{theorem}

We recall that the global attractor $\A$ of $S(t)$ acting on $\H$
is the unique compact subset of $\H$
which is at the same time fully invariant, i.e.,
$$S(t)\A=\A,\quad\forall t\geq 0,$$
and attracting, i.e.,
$$\lim_{t\to\infty}\boldsymbol{\delta}_\H(S(t)B,\A)= 0,$$
for every bounded set $B\subset\H$, where
$\boldsymbol{\delta}_\H$ denotes the standard
Hausdorff semidistance in $\H$
(see \cite{BV,HAL,TEM}).

\begin{remark}
Within our hypotheses, the regularity of $\A$ is optimal.
On the other hand, one can prove
that $\A$ is as regular as $f$ and $g$ permit.
For instance, if $f,g\in H^n$ for every $n\in\N$, then each
component of $\A$ belongs
to $H^n$ for every $n\in\N$.
\end{remark}

The proof of the theorem will be carried out
by showing a suitable (exponential) asymptotic compactness property
of the semigroup, which will be obtained exploiting a
particular decomposition
of $S(t)$
devised in \cite{GPV}. Besides, due to such a decomposition,
it is not hard to demonstrate (e.g., following
\cite{EMZ}) the existence of regular exponential attractors
for $S(t)$ having finite fractal dimension in $\H$.
As a straightforward consequence,
recalling that the global attractor is the {\it minimal} closed attracting set, we 
have

\begin{corollary}
\label{FRAC}
The fractal dimension of $\A$ in $\H$
is finite. 
\end{corollary}

\begin{remark}
In fact, having proved the existence of the absorbing set
$\BB$, we could also consider the nonautonomous case
(when $f$ and $g$ depend on time), establishing a more general
result on the existence of the global attractor for
a process of operators, provided that $f$ and $g$
fulfill suitable translation compactness properties
(see \cite{CV,CVbook} for more details). However, in that case,
the decomposition from \cite{GPV} fails to work,
and other techniques should be employed in order to establish
asymptotic compactness,
such as the $\alpha$-contraction method~\cite{HAL}
(see also \cite{EM}, where the method is applied to
a similar, albeit autonomous, problem).
\end{remark}

We now dwell on the structure of the global attractor.
To this aim, we introduce the set
$${\mathcal S}=\big\{z\in\H:S(t)z=z,\,\forall t\geq 0\big\}
$$
of stationary
points of $S(t)$, which clearly consists of
all vectors of the form $(u,0,\theta_g)$, where $u\in H^4$ is a solution to the elliptic
problem
$$
Au+\big(\beta+\|u\|^2_1\big)A^{1/2}u= f+g.
$$
The set ${\mathcal S}$ turns out to be nonempty and
$-z_g+{\mathcal S}$ is bounded in $\H^2$.
Then, the following characterization of $\A$ holds.

\begin{proposition}
\label{propSTAZ}
The global attractor $\A$
coincides with the unstable set of ${\mathcal S}$; namely,
$$\A=
\big\{z(0): z(t) \text{ is a complete trajectory of $S(t)$ and }
\lim_{t\to \infty}\|z(-t)-{\mathcal S}\|_{\H}=0\big\}.
$$
\end{proposition}

Recall that $z(t)$ is called a complete trajectory of $S(t)$
if
$$z(t+\tau)=S(t)z(\tau),\quad\forall t\geq 0,\,\forall \tau\in\R.$$

\begin{corollary}
\label{corSTAZ}
If ${\mathcal S}$ is finite, then
$$\A=
\big\{z(0):
\lim_{t\to \infty}\|z(-t)-z_1\|_{\H}
=\lim_{t\to \infty}\|z(t)-z_2\|_{\H}=0\big\},
$$
for some $z_1,z_2\in{\mathcal S}$.
If ${\mathcal S}$ consists of a single
element $z_{\rm s}\in\H^2$, then $\A=\{z_{\rm s}\}$.
\end{corollary}

As shown in \cite{CZGP}, the set ${\mathcal S}$
is always finite when all the eigenvalues $\lambda_n$ of $A$
(recall that $\lambda_n\uparrow\infty$)
satisfying the relation
$$\beta<-\sqrt{\lambda_n}\,$$
are simple
(this is the case
in the concrete problem \eqref{PROB}-\eqref{BC}), while it possesses a
single element $z_{\rm s}$ if $\beta\geq -\sqrt{\lambda_1}\,$.
In particular, we have

\begin{corollary}
\label{corSTAZ2}
If $\beta>-\sqrt{\lambda_1}\,$ and $f+g=0$, then
$\A=\{z_g\}$ and
$$\boldsymbol{\delta}_\H(S(t)B,\A)
=\sup_{z\in B}\|S(t)z-z_g\|_\H\leq\Q(\|B\|_\H)\e^{-\varkappa t},$$
for some $\varkappa>0$
and some positive increasing
function $\Q$. Both $\varkappa$ and $\Q$
can be explicitly computed.
\end{corollary}

We conclude the section discussing the injectivity
of $S(t)$ on $\A$.

\begin{proposition}
\label{BACK}
The map $S(t)_{|\A}:\A\to\A$ fulfills the backward uniqueness
property; namely, the equality $S(t)z_1=S(t)z_2$, for some $t>0$
and $z_1,z_2\in\A$, implies that $z_1=z_2$.
\end{proposition}

As a consequence, the map $S(t)_{|\A}$ is a bijection on $\A$, and so
it can be
extended to negative times by the formula
$$S(-t)_{|\A}=[S(t)_{|\A}]^{-1}.$$
In this way, $S(t)_{|\A}$, $t\in\R$,
is a strongly continuous (in the topology of $\H$) group of operators on
$\A$.

\section{Proofs of the Results}

\subsection{An equivalent problem}
Denoting as usual $S(t)z=(u(t),\pt u(t),\theta(t))$,
for some given $z=(u_0,u_1,\theta_0)\in\H$,
we introduce the function
$$\omega(t)=\theta(t)-\theta_g.$$
It then is apparent that $(u(t),\pt u(t),\omega(t))$ solves
\begin{equation}
\label{BASENEW}
\begin{cases}
\ptt u+Au-A^{1/2}\omega+
\big(\beta+\|u\|^2_1\big)A^{1/2}u= h,\\
\noalign{\vskip.7mm}
\pt \omega+A^{1/2}\omega+A^{1/2}\pt u=0,
\end{cases}
\end{equation}
where
$$h=f+g\in H,$$
with the initial conditions
$$(u(0),\pt u(0),\omega(0))=z-z_g.$$
According to Proposition~\ref{EU}, system \eqref{BASENEW} generates
a strongly continuous semigroup $S_0(t)$ on $\H$, which clearly fulfills
the relation
\begin{equation}
\label{SIM}
S(t)(\zeta+z_g)=z_g+S_0(t)\zeta,\quad\forall \zeta\in\H.
\end{equation}
Thus, from Theorem~\ref{thmABS}, we learn that $S_0(t)$
possesses
the absorbing set
$$\BB_0=-z_g+\BB.$$
Using also \eqref{TFIN}, we have the uniform bound
\begin{equation}
\label{BOUND}
\sup_{t\geq 0}\sup_{\zeta\in\BB_0}\|S_0(t)\zeta\|_\H\leq C.
\end{equation}

\begin{remark}
Here and till the end of the section, the generic constant $C$ depends
only on $\beta$, $\|h\|$ and the size of the absorbing set $\BB_0$.
\end{remark}

In light of \eqref{SIM}, Theorem~\ref{MAIN} is an immediate consequence
of the following result.

\begin{theorem}
\label{MAINNEW}
The semigroup $S_0(t)$ acting on $\H$ possesses the connected
global attractor $\A_0$ bounded in $\H^2$.
\end{theorem}

We postpone the proof of Theorem~\ref{MAINNEW}, which requires several steps.
In the sequel, for $\zeta=(u_0,u_1,\omega_0)\in\H$, we denote
$$S_0(t)\zeta=(u(t),\pt u(t),\omega(t)),$$
whose corresponding energy is given by
$$\E_0(t)=\frac12\|S_0(t)\zeta\|^2_\H+\frac14\big(\beta+\|u(t)\|_1^2\big)^2.
$$
Therefore,
the functional
$$\L_0(t)=\E_0(t)-\l h, u(t)\r
$$
satisfies the differential equality
\begin{equation}
\label{L0}
\ddt \L_0+\|\omega\|_1^2=0.
\end{equation}
It is then an easy matter to show that $\L_0$ is a Lyapunov functional
for $S_0(t)$, and by means of standard arguments (see, e.g., \cite{BV,HAL,TEM})
we conclude that
$$\A_0=
\big\{\zeta(0): \zeta(t) \text{ is a complete trajectory of $S_0(t)$ and }
\lim_{t\to \infty}\|\zeta(-t)-{\mathcal S}_0\|_{\H}=0\big\},
$$
where
$${\mathcal S}_0=-z_g+{\mathcal S}$$
is the (nonempty) set of stationary points of $S_0(t)$.
Besides, if
${\mathcal S}_0$ is finite,
$$\A_0=
\big\{\zeta(0):
\lim_{t\to \infty}\|\zeta(-t)-\zeta_1\|_{\H}
=\lim_{t\to \infty}\|\zeta(t)-\zeta_2\|_{\H}=0\big\},
$$
for some $\zeta_1,\zeta_2\in{\mathcal S}_0$.
In particular, when ${\mathcal S}_0$ consists of
a single element, recalling that the Lyapunov functional is decreasing
along the trajectories, there exists only one (constant) complete trajectory
of $S_0(t)$, so implying that
$\A_0$ is a singleton.
On account of \eqref{SIM}, this provides the proofs
of Proposition~\ref{propSTAZ} and Corollary~\ref{corSTAZ}.
In the same fashion,
Proposition~\ref{BACK} follows from the analogous statement
for $S_0(t)$, detailed in the next proposition.

\begin{proposition}
\label{BACK0}
The map $S_0(t)_{|\A_0}:\A_0\to\A_0$ fulfills the backward uniqueness
property.
\end{proposition}

\begin{proof}
We follow a classical method devised by
Ghidaglia~\cite{GHID} (see also \cite{TEM}, \S III.6), along with an
argument devised in
\cite{CHLA02}.
For $\zeta_1,\zeta_2\in\A_0$,
let us denote
$$S_0(t)\zeta_\imath=(u_\imath(t),v_\imath(t),\omega_\imath(t)),\quad \imath=1,2,$$
and
$$\zeta(t)=(u(t),v(t),\omega(t))=S_0(t)\zeta_1-S_0(t)\zeta_2.$$
Assume that $\zeta(T)=0$ for some $T>0$.
We draw the conclusion if we show that
\begin{equation}
\label{ZETA}
\zeta(0)=\zeta_1-\zeta_2=0.
\end{equation}
On account of \eqref{BASENEW}, the column vector
$$\xi(t)=(A^{1/2}u(t),v(t),\omega(t))^\top$$
satisfies the
differential equation
$$
\partial_t \xi+A^{1/2}\!\!\cdot{\mathbb B}\, \xi =\G,
$$
where we put
$${\mathbb B}=
\left(
\begin{array}{rrr}
0 & -1 &0\\
1 & 0 & -1\\
0 & 1 & 1
\end{array}
\right),\quad
\G=
\left(
\begin{matrix}
0 \\
\big(\|u_2\|_1^2-\|u_1\|_1^2\big)A^{1/2}u_2-\big(\beta+\|u_1\|_1^2\big)A^{1/2}u\\
0
\end{matrix}
\right).
$$
The matrix ${\mathbb B}$ possesses three distinct eigenvalues of strictly
positive real parts, precisely: $a\sim 0.57$
and $b\pm {\rm i} c$, with $b\sim 0.22$ and $c\sim 1.31$.
Hence, there exists a (complex) invertible $(3\times 3)$-matrix ${\mathbb U}$
such that
$${\mathbb U}^{-1}{\mathbb B}{\mathbb U}={\mathbb D},$$
where ${\mathbb D}$ is the diagonal matrix whose entries are the
eigenvalues of ${\mathbb B}$.
Accordingly, setting $\G_\star={\mathbb U}^{-1}\G$,
the (complex) function $\xi_\star(t)={\mathbb U}^{-1}\xi(t)$
fulfills
$$
\partial_t \xi_\star+A^{1/2}\!\!\cdot{\mathbb D}\, \xi_\star =\G_\star.
$$
Besides, $\xi_\star(T)=0$.
At this point, for $r=0,1$, we consider the complex Hilbert spaces
$$\W^r=H_{\mathbb C}^r\times H_{\mathbb C}^r\times H_{\mathbb C}^r,$$
$H_{\mathbb C}^r$ being the complexification of $H^r$.
It is convenient to endow $\W^1$ with the equivalent norm
$$\|\xi_\star\|_{\W^1}^2=a\|w_\star\|_{1}^2
+b\|v_\star\|_{1}^2+b\|\omega_\star\|_{1}^2,\quad
\xi_\star=(w_\star,v_\star,\omega_\star).
$$
With this choice,
$$\l A^{1/2}\!\!\cdot{\mathbb D}\, \xi_\star,\xi_\star\r_\W=\|\xi_\star\|^2_{\W^1}.
$$
It is also apparent that
$$\|\G\|_\W\leq k\|A^{1/2}u\|=k\|\xi\|_\W,$$
for some $k>0$ 
independent of $\zeta_1,\zeta_2\in\A_0$. Thus, up to redefining $k>0$,
$$\|\G_\star\|_{\W}\leq k\|\xi_\star\|_{\W}.$$
Next, we define the function
$$
\Gamma(t)=\frac{\|\xi_\star (t)\|_{\W^1}^2}{\|\xi_\star(t)\|_\W^2}.$$
Taking the time-derivative of $\Gamma$, 
and exploiting the equality
$$\l A^{1/2}\!\!\cdot{\mathbb D}\, \xi_\star-\Gamma\xi_\star,\Gamma \xi_\star\r_\W=0,
$$
we obtain 
\begin{align*}
\ddt\Gamma&=\frac{-2\| A^{1/2}\!\!\cdot{\mathbb D}\, \xi_\star
-\Gamma \xi_\star\|^2_\W}{\|\xi_\star\|_\W^2}
+\frac{2\Re\langle A^{1/2}\!\!\cdot{\mathbb D}\, \xi_\star-\Gamma \xi_\star,
\G_\star\rangle_\W}{\|\xi_\star\|_\W^2}\\
&\leq \frac{-2\| A^{1/2}\!\!\cdot{\mathbb D}\, \xi_\star
-\Gamma \xi_\star\|^2_\W}{\|\xi_\star\|_\W^2}
+\frac{\| A^{1/2}\!\!\cdot{\mathbb D}\, \xi_\star
-\Gamma \xi_\star\|^2_\W}{\|\xi_\star\|_\W^2}
+\frac{\|\G_\star\|^2_\W}{\|\xi_\star\|_\W^2}\leq k^2,
\end{align*}
and an integration in time provides the estimate 
$$
\Gamma(t)\leq \Gamma(0)+k^2t.
$$
If \eqref{ZETA} is false, by continuity, $\zeta(t)\neq 0$
in a neighborhood of zero. In turn, $\xi_\star(t)\neq 0$
in the same neighborhood.
Recalling that $\xi_\star(T)=0$, there exists $T_0\in(0,T]$
such that
$\xi_\star(t)\neq 0$ on $[0,T_0)$ and $\xi_\star(T_0)=0$.
Taking the time-derivative of $\log\|\xi_\star\|_\W^{-1}$, we 
find
$$
\ddt\log \|\xi_\star\|_\W^{-1}
=-\frac12\ddt\log \|\xi_\star\|_\W^{2}
=\Gamma -\frac{\Re\langle \xi_\star,\G_\star\rangle_\W}{\|\xi_\star\|_\W^2}
\leq \Gamma +\frac{\|\G_\star\|_\W}{\|\xi_\star\|_\W}\leq\Gamma+k.
$$
Integrating on $(0,t)$, with $t<T_0$, we conclude that
$$
\log \|\xi_\star(t)\|_\W^{-1}
\leq \log \|\xi_\star(0)\|_\W^{-1}+T_0\Gamma(0)
+k T_0+\frac12 k^2T_0^2.
$$
This
produces a uniform bound on $\log \|\xi_\star(t)\|_\W^{-1}$
over $[0,T_0)$, in contradiction with the fact
that $\|\xi_\star(T_0)\|_{\W}=0$.
\end{proof}

\subsection{Proof of Theorem \ref{MAINNEW}}
We need first to prove a suitable
dissipation integral for the norm of $\pt u$.

\begin{lemma}
\label{INTut}
For every $\nu>0$ small, there is $C_\nu>0$ such that
$$
\int_s^t \|\pt u(\tau)\|^2\d \tau\leq \nu(t-s)+C_\nu,\quad\forall t>s\geq 0,
$$
whenever $\zeta\in\BB_0$.
\end{lemma}

\begin{proof}
Let $\zeta\in\BB_0$.
An integration of \eqref{L0}, together with
\eqref{BOUND}, yield the
uniform bound
\begin{equation}
\label{omegaINTB}
\int_0^\infty\|\omega(t)\|_1^2\d t\leq C.
\end{equation}
uniformly with respect to $\zeta\in\BB_0$.
We now consider the auxiliary functionals, analogous to those
in the proof of Theorem~\ref{thmABS},
$$\Phi_0(t)=\l\pt u(t),u(t)\r,\quad
\Psi_0(t)=\l \pt u(t),\omega(t)\r_{-1},$$
which satisfy the equalities
$$
\ddt\Phi_0+\|u\|_2^2+\|u\|_1^4+\beta\|u\|_1^2
=\|\pt u\|^2+\l u,\omega\r_1+\langle u, h \rangle,
$$
and
$$\ddt\Psi_0+\|\pt u\|^2=\|\omega\|^2
-\l\pt u,\omega\r-\l u,\omega\r_1+\l\omega,h\r_{-1}
-\big(\beta+\|u\|_1^2\big)\l u,\omega\r.$$
Then, we easily get
$$
\ddt\Phi_0+\frac12 \|u\|_2^2\leq \|\pt u\|^2+C+C\|\omega\|_1^2,
$$
and, for all positive $\eps\leq 1$,
$$\ddt\Psi_0+\frac12 \|\pt u\|^2\leq
\frac{\eps}{2}\|u\|_2^2 + C \eps
+\frac{C}{\eps}\|\omega\|_1^2.$$
Therefore, for every $\eps\leq 1/4$,
$$
\ddt\big\{\eps\Phi_0+\Psi_0\big\}
+\frac14 \|\pt u\|^2\leq
C \eps+\frac{C}{\eps}\|\omega\|_1^2.$$
Integrating the last inequality on $(s,t)$,
and using \eqref{BOUND} and \eqref{omegaINTB},
we conclude that
$$\int_s^t\|\pt u(\tau)\|^2\d \tau\leq C\eps(t-s)+\frac{C}{\eps}.
$$
Setting $\nu=C\eps$ and $C_\nu=C/\eps$ the proof is completed.
\end{proof}

We shall also make use of the following Gronwall-type lemma
(see, e.g., \cite{GGMP1}).

\begin{lemma}
\label{GRNWOLD}
Let $\Lambda:\R^+\to\R^+$ be an absolutely continuous function
satisfying, for some $\nu>0$, the differential inequality
$$
\ddt\Lambda(t)+2\nu \Lambda(t)\leq \psi(t)\Lambda(t),
$$
where $\psi:\R^+\to\R^+$ is any locally
summable functions such that
$$\int_s^t \psi(\tau)\d \tau\leq \nu (t-s)+K,\quad \forall t>s\geq 0,$$
with $K\geq 0$.
Then,
$$
\Lambda(t)\leq \e^K\Lambda(0)\e^{-\nu t}.
$$
\end{lemma}

At this point, exploiting the
interpolation inequality
$$\|u\|_1^2\leq \|u\|\|u\|_2,$$
we choose $\alpha>0$ large enough such that
\begin{equation}
\label{ESTIMA}
\frac14\|u\|_2^2\leq \frac12\|u\|_2^2+\beta\|u\|_1^2+\alpha\|u\|^2\leq
M\|u\|_2^2,
\end{equation}
for some $M=M(\alpha,\beta)\geq 1$.
Then, following~\cite{GPV},
we decompose the solution $S_0(t)\zeta$, with $\zeta\in\BB_0$, into the sum
$$S_0(t)\zeta=L(t)\zeta+K(t)\zeta,$$
where
$$L(t)\zeta=(v(t),\pt v(t),\eta(t))\quad\text{and}\quad
K(t)\zeta=(w(t),\pt w(t),\rho(t))$$
are the (unique) solutions to the Cauchy problems
\begin{equation}
\label{DECAY}
\begin{cases}
\ptt v+Av-A^{1/2}\eta+
(\beta+\|u\|^2_1)A^{1/2}v+\alpha v= 0,\\
\noalign{\vskip.5mm}
\pt \eta+A^{1/2}\eta+A^{1/2}\pt v=0,\\
\noalign{\vskip.7mm}
(v(0),\pt v(0),\eta(0))=\zeta,
\end{cases}
\end{equation}
and
\begin{equation}
\label{CPT}
\begin{cases}
\ptt w+Aw-A^{1/2}\rho+
(\beta+\|u\|^2_1)A^{1/2}w-\alpha v= h,\\
\noalign{\vskip.5mm}
\pt \rho+A^{1/2}\rho+A^{1/2}\pt w=0,\\
\noalign{\vskip.7mm}
(w(0),\pt w(0),\rho(0))=0.
\end{cases}
\end{equation}
We begin to prove the exponential decay of $L(t)\zeta$.

\begin{lemma}
\label{lemmaDECAY}
There is $\varkappa>0$ such that
$$\sup_{\zeta\in\BB_0}\|L(t)\zeta\|_{\H}\leq C\e^{-\varkappa t}.$$
\end{lemma}

\begin{proof}
Denoting for simplicity
$$E_0(t)=\|L(t)\zeta\|_\H^2,$$
we define, for $\eps>0$, the functional
$$\Lambda_0(t)=\Theta_0(t)+\eps\Upsilon_0(t),$$
where
\begin{align*}
\Theta_0(t) &=E_0(t)+\beta\|v(t)\|_1^2
+\alpha\|v(t)\|^2+\|u(t)\|_1^2\|v(t)\|_1^2,\\
\Upsilon_0(t)&=\l\pt v(t),v(t)\r+2\l \pt v(t),\eta(t)\r_{-1}.
\end{align*}
It is clear from \eqref{BOUND} and \eqref{ESTIMA} that,
for all $\eps$ small enough,
\begin{equation}
\label{CTRLE0}
\frac12 E_0\leq \Lambda_0\leq C E_0.
\end{equation}
Due to \eqref{BOUND}, \eqref{DECAY} and \eqref{CTRLE0},
we have
$$
\ddt\Theta_0+2\|\eta\|^2_1=2\l A^{1/2}u,\pt u\r\|v\|_1^2\leq C\|\pt u\|\|v\|_1^2
\leq C\|\pt u\|\Lambda_0,
$$
and
\begin{align*}
&\ddt\Upsilon_0+\big(\|v\|_2^2+\beta\|v\|_1^2+\alpha\|v\|^2\big)
+\|u\|_1^2\|v\|_1^2+\|\pt v\|^2\\
\noalign{\vskip1mm}
&=2\|\eta\|^2-2\l\pt v,\eta\r-\l v,\eta\r_1
-2\big(\beta+\|u\|_1^2\big)\l v,\eta\r
-2\alpha\l v,\eta\r_{-1}\\
&\leq\frac12\|\pt v\|^2+\frac12\|v\|^2_2+C\|\eta\|^2_1.
\end{align*}
Thus, using \eqref{ESTIMA}, we obtain
$$
\ddt\Upsilon_0+\frac14\|v\|_2^2
+\frac12\|\pt v\|^2\leq C\|\eta\|^2_1.
$$
We conclude that
$$
\ddt\Lambda_0+\frac{\eps}4\|v\|_2^2
+\frac{\eps}2\|\pt v\|^2+(2-C\eps)\|\eta\|_1^2
\leq C\|\pt u\|\Lambda_0
\leq  \frac{\eps}{16}\Lambda_0+C\|\pt u\|^2\Lambda_0.
$$
Appealing again to \eqref{CTRLE0}, it is apparent that,
for all $\eps>0$ small enough, we have the inequality
$$
\ddt\Lambda_0+\frac{\eps}{16}\Lambda_0
\leq C\|\pt u\|^2\Lambda_0.
$$
The conclusion follows from Lemma~\ref{INTut}
and Lemma~\ref{GRNWOLD}, using again
\eqref{CTRLE0}.
\end{proof}

\begin{corollary}
\label{corDIRDEC}
If $\beta>-\sqrt{\lambda_1}\,$ and $h=0$,
then
$$\|S_0(t)\zeta\|_\H\leq\Q(\|\zeta\|_\H)\e^{-\varkappa t},$$
for some $\varkappa>0$
and some positive increasing
function $\Q$. Both $\varkappa$ and $\Q$
can be explicitly computed.
\end{corollary}

\begin{proof}
We preliminarily observe that it suffices to prove the result
for $\zeta\in\BB_0$. Moreover,
choosing a constant $\sigma$ such that
$$-\beta\lambda_1^{-1/2}<\sigma<1,$$
we find the controls
\begin{equation}
\label{ESTIMDEC}
m\|u\|_2^2
\leq\sigma \|u\|_2^2+\beta\|u\|_1^2 \leq
\big(1+|\beta|\lambda_1^{-1/2}\big)\|u\|_2^2,
\end{equation}
with
$$m=
\begin{cases}
\sigma&\text {if }\beta\geq 0,\\
\sigma+\beta\lambda_1^{-1/2}&\text {if }\beta<0.
\end{cases}
$$
Then, we just recast the proof of Lemma~\ref{lemmaDECAY},
with \eqref{ESTIMDEC} in place of \eqref{ESTIMA}.
\end{proof}

\begin{remark}
Exploiting \eqref{SIM}, from Corollary~\ref{corDIRDEC}
we readily get the proof of Corollary~\ref{corSTAZ2}.
\end{remark}

The next lemma shows the uniform boundedness of $K(t)\BB_0$ in
the more regular space $\H^2$, compactly embedded into $\H$.

\begin{lemma}
\label{lemmaCPT}
The estimate
$$\sup_{\zeta\in\BB_0}\|K(t)\zeta\|_{\H^2}\leq C$$
holds for every $t\geq 0$.
\end{lemma}

\begin{proof}
We first observe that, from \eqref{BOUND} and Lemma~\ref{lemmaDECAY},
we have
\begin{equation}
\label{ienki}
\|w\|_3^2\leq \|w\|_2\|w\|_4
\leq (\|u\|_2+\|v\|_2)\|w\|_4\leq C\|w\|_4.
\end{equation}
Setting
$$E_1(t)=\|K(t)\zeta\|_{\H^2}^2,$$
we define, for $\eps>0$, the functional
$$\Lambda_1(t)=\Theta_1(t)+\eps\Upsilon_1(t),$$
where
\begin{align*}
\Theta_1(t) &=E_1(t)+\big(\beta+\|u(t)\|_1^2\big)\|w(t)\|_3^2-2\l Aw(t),h\r,\\
\Upsilon_1(t)&=\l\pt w(t),w(t)\r_2+2\l \pt w(t),\rho(t)\r_{1}.
\end{align*}
Note that, from \eqref{BOUND} and \eqref{ienki},
\begin{equation}
\label{CTRLE1}
\frac12 E_1-C\leq \Lambda_1\leq C E_1+C,
\end{equation}
for all $\eps$ small enough.
Exploiting Lemma~\ref{lemmaDECAY},
\eqref{BOUND} and \eqref{ienki},
we have
$$
\ddt\Theta_1+2\|\rho\|^2_3=2\alpha\l v,\pt w\r_2+2\l A^{1/2}u,\pt u\r\|w\|_3^2
\leq\frac{\eps}{4}\|w\|_4^2+
\frac{\eps}{4}\|\pt w\|_2^2+\frac{C}{\eps},
$$
and
\begin{align*}
&\ddt\Upsilon_1+\|w\|_4^2+\|\pt w\|_2^2
+\|u\|_1^2\|w\|_3^2
\\
\noalign{\vskip1mm}
&=-\beta\|w\|_3^2+2\|\rho\|_2^2
-\l w,\rho\r_3
-2\big(\beta+\|u\|_1^2\big)\l w,\rho\r_2-2 \langle \pt w,\rho \rangle_2\\
&\quad +\alpha\l v,w\r_2+2\alpha\l v,\rho\r_1+\l Aw+2 A^{1/2}\rho,h\r\\
&\leq \frac14\|w\|_4^2+\frac14\|\pt w\|_2^2+C\|\rho\|_3^2+C,
\end{align*}
which entails
$$
\ddt\Upsilon_1+\frac34 \|w\|_4^2+\frac34\|\pt w\|_2^2
\leq C\|\rho\|_3^2+C.
$$
Collecting the above estimates, we are led to
$$
\ddt\Lambda_1+\frac{\eps}2 \|w\|_4^2+\frac{\eps}{2}\|\pt w\|_2^2
+(2-C\eps)\|\rho\|_3^2\leq \frac{C}{\eps}.
$$
On account of \eqref{CTRLE1}, we can now fix $\eps$ small enough
such that the inequality
$$
\ddt\Lambda_1+\nu\Lambda_1\leq C
$$
holds for some $\nu>0$.
Applying the Gronwall lemma, and using again \eqref{CTRLE1}, we are done.
\end{proof}

In conclusion, Lemma~\ref{lemmaDECAY} and Lemma~\ref{lemmaCPT}
tell us that $S_0(t)\BB_0$ is (exponentially) attracted
by a bounded subset of $\H^2$, thus, precompact in $\H$.
As is well known from the theory of dynamical systems
(see, e.g., \cite{BV,HAL,TEM}), this yields the existence
of the global attractor $\A_0$, bounded in $\H^2$, for the semigroup
$S_0(t)$ acting on $\H$.
The proof of Theorem~\ref{MAINNEW} is finished.

\section{A More General Model}

\noindent
In this final section, we discuss a more general abstract
problem, obtained by adding the term
$\gamma A^{1/2}\ptt u$ to the first equation of \eqref{BASE}, where
$\gamma\geq 0$ is the so-called rotational parameter.

Given $\gamma\geq 0$,
we define the strictly positive selfadjoint operator on $H$
$$M_\gamma=1+\gamma A^{1/2},$$
with domain 
$\D(M_\gamma)=H^2$ (when $\gamma>0$). 
Since the operator $M_\gamma$ commutes with $A$
and all its powers, we introduce the spaces
$$H^r_\gamma=\D(A^{(r-1)/4}M_\gamma^{1/2}),\quad r\in\R,$$
with inner products and norms
$$\l u,v\r_{r,\gamma}=\l A^{(r-1)/4}M_\gamma^{1/2}u,A^{(r-1)/4}M_\gamma^{1/2}v\r,\quad
\|u\|_{r,\gamma}=\|A^{(r-1)/4}M_\gamma^{1/2}u\|.$$
Finally, we set
$$\V^r_\gamma=H^{r+2}\times H^{r+1}_\gamma\times H^r.
$$
Again, we agree to omit the index $r$ when $r=0$.

\begin{remark}
Note that
\begin{equation}
\label{DUENORME}
\|u\|^2_{r,\gamma}=\|u\|^2_{r-1}+\gamma\|u\|^2_{r}
\leq \Big(\frac1{\sqrt{\lambda_1}}\,+\gamma\Big)\|u\|^2_{r}.
\end{equation}
Hence,
when $\gamma>0$, the space $H^{r}_\gamma$ is just $H^r$
endowed with an equivalent norm, whereas
$H^r_0=H^{r-1}$ and 
$\V^r_0=\H^r$. 
\end{remark}

We consider the Cauchy problem on $\V_\gamma$
\begin{equation}
\label{ROTSYS}
\begin{cases}
M_\gamma\ptt u+Au-A^{1/2}\theta+
\big(\beta+\|u\|^2_1\big)A^{1/2}u= f(t),\quad t>0,\\
\noalign{\vskip.7mm}
\pt \theta+A^{1/2}\theta+A^{1/2}\pt u=g(t),\quad t>0,\\
u(0)=u_0,\quad
\pt u(0)=u_1,
\quad\theta(0)=\theta_0,
\end{cases}
\end{equation}
of which \eqref{BASE} is just the particular instance corresponding to $\gamma=0$.
In concrete models, the additional term $\gamma  A^{1/2} \ptt u$ accounts for
the presence of rotational inertia.
With $f$ and $g$ as in Proposition~\ref{EU}, this system
generates a family of solution operators
$S^\gamma(t)$ on $\V_\gamma$, satisfying the joint
continuity property
$$(t,z)\mapsto S^\gamma(t)z\in\C(\R^+\times\V_\gamma,\V_\gamma).$$
The energy at time $t$ corresponding
to the initial data $z\in\V_\gamma$ now reads
$$\E^\gamma(t)=\frac12\|S^\gamma(t)z\|^2_{\V_\gamma}+\frac14\big(\beta+\|u(t)\|_1^2\big)^2,
$$
and the energy identity \eqref{E} is still true replacing $\E$ with $\E^\gamma$.
As a matter of fact, all the results stated in the previous sections
extend to the present case.

\begin{theorem}
\label{ALL}
Theorems~\ref{thmABS}, \ref{MAIN}, Corollaries~\ref{FRAC}, \ref{corSTAZ}, \ref{corSTAZ2}
and Proposition~\ref{propSTAZ} continue to hold with $S^\gamma(t)$ 
and $\V^r_\gamma$ in place of $S(t)$ and $\H^r$.
\end{theorem}

\begin{proof}[Sketch of the proof]
Repeat exactly the same demonstrations,
simply replacing (clearly, besides $S(t)$ with $S^\gamma(t)$
and $\H^r$ with $\V^r_\gamma$)
each occurrence
of $\pt u$ [or $\pt v$, $\pt w$] with $M_\gamma\pt u$
[or $M_\gamma\pt v$, $M_\gamma\pt w$]
in the definitions of the auxiliary functionals
$\Phi$, $\Psi$, $\Phi_0$, $\Psi_0$, $\Upsilon_0$, $\Upsilon_1$.
The integral estimate of Lemma~\ref{INTut}
improves to 
$$
\int_s^t \|\pt u(\tau)\|^2_{1,\gamma}\d \tau\leq \nu(t-s)+C_\nu,\quad\forall t>s\geq 0,
$$
although, as we will see in a while, the original
estimate would suffice.
For example, let us examine
more closely the modifications needed in the new proofs
of Lemma~\ref{lemmaDECAY} and Lemma~\ref{lemmaCPT}. We keep the same notation, 
just recalling that now the terms $\ptt v$ and
$\ptt w$ in \eqref{DECAY}-\eqref{CPT}
are replaced by $M_\gamma\ptt v$ and $M_\gamma\ptt w$, respectively.
The estimate on
$\frac{\d}{\d t}\Theta_0$ remains unchanged. Concerning 
$\frac{\d}{\d t}\Upsilon_0$, the term $\|\pt v\|^2$ in the left-hand 
side turns into $\|\pt v\|_{1,\gamma}^2$, and in the right-hand side
we have $-2\l M_\gamma \pt v,\eta\r$ instead of $-2\l \pt v,\eta\r$.
But thanks to \eqref{DUENORME},
$$-2\l M_\gamma \pt v,\eta\r
\leq 2\|\pt v\|_{1,\gamma}\|\eta\|_{1,\gamma}
\leq \frac12\|\pt v\|_{1,\gamma}^2+C\|\eta\|_1^2,$$
for some $C>0$ independent of $\gamma$, provided that, say, $\gamma\leq 1$.
So, we end up with the same differential inequality for $\Lambda_0$,
which yields the desired claim in light of the dissipation integral
for $\|\pt u\|$.
Coming to Lemma~\ref{lemmaCPT}, 
we readily have (cf.\ \eqref{DUENORME}) 
$$
\ddt\Theta_1+2\|\rho\|^2_3
\leq\frac{\eps}{4}\|w\|_4^2+
\frac{\eps}{4}\|\pt w\|_2^2+\frac{C}{\eps}
\leq\frac{\eps}{4}\|w\|_4^2+
\frac{\eps}{4}\|\pt w\|_{3,\gamma}^2+\frac{C}{\eps},
$$
whereas in the estimate of $\frac{\d}{\d t}\Upsilon_1$
the term $\|\pt w\|^2_2$ in the left-hand 
side becomes $\|\pt w\|_{3,\gamma}^2$, and in the right-hand side
$-2\l M_\gamma \pt w,\rho\r_2$ substitutes $-2\l \pt w,\rho\r_2$.
A further use of \eqref{DUENORME} entails the control
$$-2\l M_\gamma \pt w,\rho\r_2
\leq 2\|\pt w\|_{3,\gamma}\|\rho\|_{3,\gamma}
\leq \frac14\|\pt w\|_{3,\gamma}^2+C\|\rho\|_3^2.$$
Once again, we are led to the same differential inequality for $\Lambda_1$.
\end{proof}

In particular, when $\gamma>0$, the global attractor $\A^\gamma$ of the semigroup
$S^\gamma(t)$ is a bounded subset of 
$\V_\gamma^2=H^4\times H^3\times H^2$. This improves the conclusions of
\cite{BUC} where, for $\gamma>0$, the boundedness of $\A^\gamma$ is obtained 
only in the 
{\it intermediate} space $H^3\times H^2\times H^2$. A straightforward albeit relevant consequence
of the $\V^2_\gamma$-regularity is emphasized in the next corollary.

\begin{corollary}
For every $\gamma\geq 0$, the solutions to \eqref{ROTSYS}
with initial data 
on the attractor are strong solutions, i.e., the equations
hold almost everywhere.
\end{corollary}

It is worth noting that
all the estimates obtained in the
proofs are {\it uniform} with respect to $\gamma$
(assuming $\gamma$ bounded from above).
Indeed, the dependence on the rotational parameter enters only
through the definition of the norm.
Then, recasting a standard argument from \cite{HAL},
the family $\{\A^\gamma\}$ is easily shown to
be upper semicontinuous at $\gamma=0$, namely,
$$\lim_{\gamma\to 0}\boldsymbol{\delta}_\H(\A^\gamma,\A)= 0.$$
On the contrary, the analogue of Proposition~\ref{BACK} does not seem to follow
by a straightforward adaptation of the preceding argument. However, the backward uniqueness
property on the attractor holds for $\gamma>0$
as well, and it can be proved
as in \cite{CHLA08}.

\subsection*{Acknowledgments}
We are grateful to the Referee for several valuable suggestions and comments.



\end{document}